\documentclass{amsart}

\usepackage{amssymb}

\usepackage{graphicx}

\usepackage[cmtip,all]{xy}

\newtheorem{theorem}{Theorem}[section]
\newtheorem{lemma}[theorem]{Lemma}

\newcommand{\tuborg}{\left\{\begin{array}{ll}}
\newcommand{\sluttuborg}{\end{array}\right.}
\newcommand{\calO}{\mathcal{O}}

\newcommand{\calZ}{\mathcal{Z}}

\newcommand{\calM}{\mathcal{M}}

\newcommand{\calB}{\mathcal{B}}
\newcommand{\calC}{\mathcal{C}}

\newcommand{\bbZ}{\mathbb{Z}}
\newcommand{\bbP}{\mathbb{P}}
\newcommand{\bbQ}{\mathbb{Q}}

\newcommand{\bbC}{\mathbb{C}}

\newcommand{\bbA}{\mathbb{A}}

\newcommand{\bbW}{\mathbb{W}}
\newcommand{\ord}{{\rm ord}}

\newcommand{\supp}{{\rm Supp}}

\newcommand{\spec}{{\rm Spec}}

\newcommand{\sgn}{{\rm sgn}}
\newcommand{\res}{{\rm res}}
\newcommand{\tr}{{\rm Tr}}

\theoremstyle{definition}
\newtheorem{definition}[theorem]{Definition}
\newtheorem{example}[theorem]{Example}

\newtheorem{corollary}[theorem]{Corollary}
\newtheorem{conj}[theorem]{Conjecture}
\newtheorem{question}[theorem]{Question}

\theoremstyle{remark}
\newtheorem{remark}[theorem]{Remark}

\numberwithin{equation}{section}

\begin{document}

\title{Regulators on additive higher Chow groups}

\author{Jinhyun Park}
\address{Department of Mathematics, Purdue University, 150 North University Street, West Lafayette, Indiana 47907, USA}
\email{jinhyun@math.purdue.edu}

\subjclass[2000]{Primary 14C15; Secondary 19D55}

\date{December 13, 2006; revised March 2, 2008}

\commby{Karen Smith}

\begin{abstract}As an attempt to understand motives over $k[x]/(x^m)$, we define the cubical additive higher Chow groups with modulus for all dimensions extending the works of S. Bloch, H. Esnault and K. R\"ulling on $0$-dimensional cycles. We give an explicit construction of regulator maps on the groups of $1$-cycles with an aid of the residue theory of A. Parshin and V. Lomadze.
\end{abstract}

\maketitle

\section*{Introduction}
This paper considers the problem of defining the motivic cohomology groups of non-reduced schemes over a field of characteristic zero. One would hope that for any noetherian scheme $X$ there may be the motivic cohomology groups $H_{\calM} ^{*} (X; \bbZ(*))$ and a spectral sequence
\begin{equation}\label{spectral sequence}E^2 _{p,q}=H_{\calM} ^{q-p} (X;\bbZ(q)) \Rightarrow K_{p+q} (X)\end{equation} that converges to the algebraic $K$-theory of D. Quillen. Indeed, when $X$ is a smooth variety over a field $k$, the higher Chow groups $CH^i (X, n)$ of S. Bloch play the role of the motivic cohomology groups (\cite{B1,V}) via $H_{\calM} ^{2i-n} (X; \bbZ (i)) = CH^i (X, n).$ But this identification fails to stay correct for non-reduced or singular schemes. For instance, when $X$ is a fat point $\spec ( k[\epsilon] ) = \spec ( k[x]/(x^2))$, the map $X_{\rm{red}} = \spec (k) \hookrightarrow X$ induces an isomorphism of $CH^i (X_{\rm{red}}, n)$ with $CH^i (X, n)$, whereas for the $K$-groups we have a splitting
\begin{equation*} K_n (k[\epsilon]) \simeq K_n (k) \oplus K_n (k[\epsilon], (\epsilon)),\end{equation*} thus the $K$-groups of $k$ and $k[\epsilon]$ differ by some relative $K$-groups. These values are known from works of many mathematicians including W. van der Kallen, F. Keune, J. Stienstra, and T. Goodwillie (\cite{Good,Keune,Sti,VdK,VdK2}), to name a few. The following important result, that guides the development of this paper, is proved by L. Hesselholt as Theorem 10 in \cite{He}:

\begin{theorem}Suppose that $A$ is a regular noetherian ring containing $\bbQ$. Let $n \geq 1$ be an integer. Then there is an isomorphism
\begin{equation}\label{Hesselholt}
K_{n+1} (A[x]/(x^m), (x)) \simeq \bigoplus_{p \geq 0} ( \Omega_{A/\bbZ} ^{n-2p} )^{m-1}.
\end{equation}
\end{theorem}

The goal of this manuscript is to show that the method of cubical additive Chow theory, initiated by S. Bloch and H. Esnault in \cite{BE2} on $0$-dimensional cycles, has a higher dimensional generalization that complements this deficiency of the \emph{usual} higher Chow groups for non-reduced schemes $\spec (k[x]/(x^m) )$. This generalization should supply the missing relative part of the motivic cohomology groups that relate to the known values of the relative $K$-groups. For instance, when $m=2$ and ${\rm char} (k) =0$, the spectral sequence \eqref{spectral sequence} and the isomorphism \eqref{Hesselholt} suggest an identification
\begin{equation*}
H^{n-2p+1} _{\calM} (k[\epsilon], (\epsilon);\bbZ (n-p+1)) \overset{?}{\simeq} \Omega_{k/\bbZ} ^{n-2p},
\end{equation*} and when $p=0$, in this paper's notation, the additive Chow group $ACH_0 (k, n; 2)$ of zero cycles (see Definitions in \S \ref{section:2}) serves the role of the group on the left as proven by S. Bloch and H. Esnault in \cite{BE2}.  More general cases of $m \geq 2$ were considered by K. R\"ulling in \cite{R} where he proves that $ACH_0 (k, n; m) \simeq \bbW_{m-1} \Omega_k ^{n}.$ The right hand side is the generalized de Rham-Witt forms of Bloch-Deligne-Hesselholt-Illusie-Madsen. 

This paper attempts the case $p=1$ using one dimensional cycles, and to justify the validity of our definition, we need a result like the above. The main result of the paper lies on this line:

\begin{theorem} Let $k$ be a field of characteristic $0$. Then there are regulator maps $$R_{2,m}: ACH_1 (k, n; m) \to \Omega_{k/\bbZ} ^{n-2}.$$
\end{theorem}
This theorem follows from the Theorem \ref{regulator1}. We use the word \emph{regulator} in its broadest sense that includes Chern character maps, cycle class maps, realizations, etc. These regulator maps can be seen as the additive analogues of the Beilinson regulator maps from the motivic cohomology groups of the smooth varieties over a field $k$ contained in $\bbC$ to the real Deligne cohomology groups (see for example \cite{G1, G4,Kerr,Kerr2}).

When $n=m=2$, the map $R_{2,2}$ can be regarded as the additive version of the Bloch-Wigner function (\cite{B3}). It is also compatible with the regulator $\rho$ of S. Bloch and H. Esnault in \cite{BE2} (see Remark \ref{reason}). 

The particular piece $R_{2,2}$ has some further applications discussed in detail in \cite{P2}; first, we use it to prove that the group $ACH_1 (k, 2;2)$ is nontrivial. Conjecturally we believe that it is isomorphic to $k$ and the regulator map is an isomorphism. Secondly, under this assumption, for an algebraically closed field $k$ of characteristic zero we can construct the additive $4$-term motivic sequence
$$
\xymatrix{ 0 \ar[r] & k \ar[r]  & T \calB_2 (k) \ar[r] & k \otimes_{\bbZ} k^{\times} \ar[r] & \Omega_{k/\bbZ} ^1 \ar[r] & 0 }$$ whose $K$-theoretic version appeared in \cite{BE2}. See also \cite{C,G3} for related discussions on dilogarithm and Euclidean scissors congruence.

\subsection*{Brief description of each section} In \S \ref{section:1}, we recall a theory of residues. The regulators of the main theorem depend on this theory. This extended notion of residues works with forms on singular varieties with higher order poles as well. Yet the reciprocity theorem, that the sum of residues over a complete variety is zero, is still valid. The most essential definitions and results are included for the convenience of the reader.

In \S \ref{section:2}, we define the cubical \emph{additive} higher Chow complex. Here is the idea: for the \emph{usual} cubical higher Chow complex (\cite{B2,T}) we consider the ambient space $X \times \square ^{n+1}$ with $\square = ( \mathbb{P} ^1 \backslash \{ 1 \} , \{ 0 , \infty \} )$, where $\{ 0, \infty \}$ are the faces. The admissible cycles are those intersecting all lower dimensional faces properly. For the cubical {\it additive} higher Chow complex, we first replace the space $\square ^{n+1}$ by the isomorphic space $\square^1 _t \times \square ^{n} = ( \bbA^1 , \{ 0, t \} ) \times \square^{n}$ for any $t \not = 0$. Then imagine the situation where $t \to 0$ so that two faces $\{0, t \}$ are collapsed into the face $2(0)$ with multiplicity $2$. We thus obtain the space $\Diamond_{n}:= ( \bbA^1 , 2 (0) ) \times \square ^{n}$. More generally, we can consider the space $\Diamond_{n}= ( \bbA^1 , m (0) ) \times \square ^{n}$ with $m \geq 2$. In addition to the proper intersection, we also require the \emph{modulus condition} (see Definition \ref{cycle groups} -(B)-(3)) that involves some data coming from intersections with the degenerate face $m (0)$. The boundary maps $\partial$ are defined in terms of the faces of the remaining coordinates on $\square ^{n}$.

In \S \ref{section:3}, when $k$ is a field of characteristic $0$, we give an explicit construction of regulators $R_{2,m}$ on the group of admissible $1$-cycles in the $(n+1)$-dimensional space $\Diamond_n$. They are sums of the residues of certain absolute K\"ahler differential $(n-1)$-forms. The main theorem comes from the reciprocity of residues. Some explicit calculations of the regulators are done as an example.

\bigskip

\noindent \textbf{Acknowledgement} This paper is based on some chapters of my doctoral thesis at the University of Chicago. I wish to thank H\'el\`ene Esnault, Stefan M\"uller-Stach, and Eckart Viehweg for their hospitality during my visit to the Universit\"at Duisburg-Essen in 2003, and Byungheup Jun and JongHae Keum for supporting my stay in 2005 at KIAS. I thank Donu Arapura, Masanori Asakura, Alexander Beilinson, Alexander Goncharov, Madhav Nori, and Tomohide Terasoma for their comments, and Matt Kerr and Kay R\"ulling for copies of their theses. I appreciate Jonathan Wahl's kind suggestion. Wilberd van der Kallen pointed out an error, and the editor and the referees of the American Journal of Mathematics offered invaluable helps in improvement of this paper. Most importantly, I wish to mention cordial supports and encouragement I received from Spencer Bloch and I thank him for being a good and patient advisor.

\section{Residues of Parshin and Lomadze}\label{section:1}

The heart of the paper is the construction of the regulator maps in \S \ref{section:3} in terms of residues of certain rational absolute K\"ahler differential forms. We face two technical challenges; we need a notion of residues for singular varieties and the notion should also work with forms with higher order poles. Nevertheless, we still want the reciprocity for complete varieties.
 
 The existence of such residues was well-known for curves, and for higher dimensional varieties it was done by El Zein, Beilinson, Parshin, and Lomadze, etc. (\cite{Be,E,L3,L4,P3}). It was improved by Yekutieli in \cite{Y}. We included some relevant definitions and results from \cite{Y} without details. In this section, $k$ is a perfect field.

\subsection{Residues at a nonsingular closed point}

Let $X$ be a variety of dimension $n$ over $k$ and let $p \in X$ be a nonsingular closed point. As $\calO_{X,p}$ is a regular local ring of dimension $n$ with the residue field $k(p)$, the Cohen structure theorem says its completion $\widehat{\calO}_{X,p}$ is isomorphic to $k(p)[[t_1, \cdots, t_n]]$. Consider the field inclusion $k(X) \hookrightarrow k(p)((t_1, \cdots, t_n))=:E$. The residue map $\res_p$ is the composition
 
 $$\res_p: \Omega_{k(X)/k} ^n = k(X) dt_1 \wedge \cdots \wedge dt_n \hookrightarrow E dt_1 \wedge \cdots \wedge dt_n \to k(p) \overset{\tr_{k(p)/k}}{\longrightarrow} k,$$
where the second arrow is defined to be
 \begin{eqnarray*}\sum_{(i_1, \cdots, i_n)\in \bbZ^n} a(i_1, \cdots, i_n) t_1 ^{i_1} \cdots t_n ^{i_n} dt_1 \wedge \cdots \wedge dt_n \mapsto  a( -1, \cdots -1)\end{eqnarray*}

It does not depend on the choice of local coordinates $ t_1, \cdots, t_n$ up to sign. 

\subsection{Parshin-Lomadze residue and reciprocity}

When $p$ lies on the singular locus, we replace the point by a sequence of specializations of points of $X$ called a chain. This approach works when the point $p$ is not closed, provided that we choose a pseudo-coefficient field of $p$ defined below. (For points $x, y \in X$, the notation $x >y$ means that $y$ is a specialization of $x$.)

\begin{definition}[\emph{c.f.} Def. 3.1.1 in \cite{Y}]Let $X$ be of finite type over $k$. A \emph{chain of length $r$} is a sequence $\xi=(x_0 , \cdots , x_r)$ of strict specializations $x_i > x_{i+1}$ of points of $X$. It is \emph{saturated} if the codimension of $\{ x_{i+1} \} ^-$ in $\{ x_i \} ^-$ is $1$ for each $i$.

 Given two chains $\xi = (x_0, \cdots, x_r)$ and $\eta= (y_0, \cdots, y_s)$ with $x_r > y_0$, their \emph{concatenation} $\xi \vee \eta$ is the chain $ (x_0, \cdots, x_r, y_0, \cdots, y_s)$.
\end{definition}

\begin{definition}[\emph{c.f.} Def. 4.1.1 in \cite{Y}]\label{coefficient field}Let $(A, \mathfrak{m})$ be a local $k$-algebra. A \emph{pseudo-coefficient field} (resp. \emph{coefficient field}) for $A$ is a $k$-algebra homomorphism $\sigma: K \to A$ where $K$ is a field and the field extension $\overline{\sigma}: K \to A \to A/\mathfrak{m}$ is finite (resp. bijective).

If $A = \widehat{\calO}_{X, y}$ for some point $y \in X$, we say that $\sigma$ is a pseudo-coefficient field (resp. coefficient field) for $y$. For a chain $\xi = (\cdots, y)$ that ends with $y$, a pseudo-coefficient field $\sigma$ for $y$ is also said to be a pseudo-coefficient field for $\xi$.
\end{definition}

If $y$ is a closed point of $X$, then there exists a unique pseudo-coefficient field from $k$. It is not usually unique for non-closed points.

\bigskip

A detailed description of this theory is not the point of this paper. For a saturated chain $\xi= (x_0, \cdots, x_r)$ of length $r$ and a pseudo-coefficient field $\sigma: K \to \widehat{\calO}_{X,y}$, we just comment that the crucial idea of the construction of the residue $\res_{\xi, \sigma} : \Omega_{k(x)/ k} ^* \to \Omega_{K/k} ^{*-r}$ is to regard the point $x_i$ as a curve over its immediate specialization $x_{i+1}$, and apply successive normalizations of local domains $\calO_{x_i, x_{i+1}}$ of dimension $1$. Interested readers should consult \cite{L4,Y} for details.

\bigskip

\noindent \textbf{Notational conventions} 

\noindent $\bullet$ When $p \in X$ is closed, consider the pseudo-coefficient field from $k$. For a chain $\xi$ ending with a closed point, we use $\res_{\xi, k}$. We sometimes write $\res_{\xi}$ without specifying the coefficient field.

\noindent $\bullet$ When $X$ is an irreducible curve over $k$ and $p$ is its closed point, we may write $\res_p$ instead of $\res_{(X, p)}$.

\noindent $\bullet$ For a chain $\xi = (x_0, \cdots, x_n)$ on $X$ with $X_i = \{ x_i \} ^-$, we may write $\res_{(X_0, X_1, \cdots, X_n)}$ instead of $\res_{\xi}$. Keep in mind that the residue is defined locally near $x_n$.

\bigskip

A rational differential form $\alpha \in \Omega^*_{k(x)/k}$ is holomorphic (Def. 4.2.3, Prop. 4.2.14 in \cite{Y}) along all but finitely many saturated chains $\xi = (x, \cdots)$ beginning with $x$. If $\alpha$ is holomorphic along $\xi$, then for any pseudo-coefficient field $\sigma$ for $\xi$ the residue $\res_{\xi, \sigma} (\alpha)$ vanishes. Thus, the sum of residues of a fixed form over a set of chains makes sense. The reciprocity theorem for residues is the following:
 
\begin{theorem}[Thm. 4.2.15 in \cite{Y}; \emph{c.f.} Thm. 3 in \cite{L4}]\label{Parshin-Lomadze} Let $X$ be a scheme of finite type over a perfect field $k$.
\begin{enumerate}
\item [(1)] Let $\xi = (\cdots, x)$ and $\eta = (y, \cdots, z)$ be saturated chains in $X$ such that $x>y$ and the codimension of $\{y\}^-$ in $\{x\}^-$ is $2$. Let $\sigma: K \to \widehat{\calO}_{X, z}$ be a pseudo-coefficient field. Then,
\begin{eqnarray*}\sum_{w \in X, x > \omega > y} \res_{\xi \vee (w) \vee \eta, \sigma} = 0.\end{eqnarray*}
\item [(2)] Suppose that $X$ is proper over $k$, and let $\xi = (\cdots, x)$ be a saturated chain on $X$ with $\dim \{ x \} ^- = 1$. Then,
\begin{eqnarray*}\sum_{w \in X, x > w } \res_{\xi \vee (w),k} = 0.\end{eqnarray*}
\end{enumerate}
\end{theorem}

Note that $X$ doesn't have to be nonsingular. Properness is essential for (2), but not for (1). Both (1) and (2) are used in the proof of the main theorem, Theorem \ref{regulator1}.
 
 \begin{remark}\label{order of pole}For a rational form $\alpha \in \Omega^* _{k(x)/k}$, there is also a notion of the \emph{order of the pole} (Def. 4.2.10 in \cite{Y}) along a saturated chain $\xi = (x, \cdots)$.
 
 When the order is $\leq 1$, we say that $\alpha$ has a \emph{simple pole} along $\xi$. In this case, the residue $\res_{\xi, \sigma} (\alpha)$ is independent of the choice of the pseudo-coefficient field $\sigma$ for $\xi$ (Cor. 4.2.13 in \cite{Y}) even when $\xi$ ends with a non-closed point. Thus, we simply write $\res_{\xi} (\alpha)$ dropping $\sigma$. This residue is called the \emph{Poincar\'e residue of $\alpha$ along $\xi$}.
 
  \end{remark}
 
 Let us state one more theorem used in this paper:
 
 \begin{theorem}[\emph{c.f.} Thm. 2 on p.515 in \cite{L4}]\label{residue under finite morphism} Let $X, Y$ be integral schemes of finite type over $k$ of the same dimension. Let $\phi: Y \to X$ be a surjective proper $k$-morphism and let $K = k(X)$, $L = k(Y)$. Suppose that $[L:K] < \infty$. Then, for any chain $\xi$ on $X$, 
 \begin{equation*}
 \res_{\xi} \circ_{\tr_{L/K}}= \sum_{\eta \to \xi} \res_{\eta},
 \end{equation*}
 where the sum is over all chains $\eta$ on $Y$ mapped to the chain $\xi$ point-wise.
 \end{theorem}
 
 \subsection{Residues on absolute K\"ahler differentials} The previous discussion generalizes to the absolute K\"ahler differentials as follows. (\emph{c.f.} Prop. 6.2 in \cite{BE2}) Let $X$ be an irreducible variety of dimension $d$ over $k$. When $r \geq d$, the natural map $\Omega_{k(X)/\bbZ} ^r \to \Omega_{k/\bbZ} ^{r-d} \otimes_k \Omega_{k(X)/k} ^d$ yields
 \begin{equation*}
 \phi: \Gamma(X, \Omega_{k(X)/\bbZ} ^r) \to \Gamma(X, \Omega_{k/\bbZ} ^{r-d} \otimes_k \Omega_{k(X)/k} ^d ) = \Omega_{k/\bbZ} ^{r-d} \otimes_k \Gamma(X, \Omega_{k(X)/ k} ^d).
 \end{equation*} Define ``$\res$'' on the absolute K\"ahler differetials as $(1 \otimes \res) \circ \phi$, where ``$\res$'' in $1 \otimes \res$ is the previously discussed residue map. When $r <d$, the natural map $\bbZ \to k(X)$ factors as $\bbZ \to k \to k(X)$, and it induces a natural map $\Omega_{k(X)/ \bbZ} ^r \to \Omega_{k(X)/k} ^r$. It yields
\begin{equation*} \phi:\Gamma(X, \Omega_{k(X)/\bbZ} ^r) \to \Gamma(X, \Omega_{k(X)/ k} ^r),\end{equation*} and we define ``$\res$'' similarly. The reciprocity for the (old) residue map immediately implies the reciprocity for this residue map on the absolute K\"ahler differentials.

\section{Cubical additive higher Chow complex}\label{section:2}

\subsection{Weil divisors on normal varieties} Let $X$ be a \emph{normal} irreducible variety over a field $k$. A \emph{Weil divisor} is an element of the free abelian group ${\rm Div} (X)$ generated by the set ${\rm PDiv}(X)$ of \emph{prime Weil divisors}, \emph{i.e.}  codimension $1$ irreducible subvarieties of $X$. If $D =\sum_Y a_Y Y$, where $ a_Y \in \bbZ$, is a Weil divisor, then the assignment $D \mapsto a_Y$ gives a homomorphism $\ord_Y: {\rm Div} (X) \to \bbZ$. We have another, but related order function $\ord_Y : k(X) \to \bbZ$ along $Y$ on the function field of $X$ as $\calO_{X,Y}$ is a DVR, $X$ being normal. A Weil divisor $D$ is called \emph{effective} if $\ord_Y D \geq 0$ for all $Y\in {\rm PDiv}(X)$, and we write $D \geq 0$.

\begin{definition}\label{weil}{\rm Let $X$ be an irreducible normal variety over a $k$ and let $D$ be a Weil divisor on $X$.

\begin{enumerate}
\item The \emph{support} of $D$, denoted by $\supp (D)$, is the set of all prime Weil divisors $Y$ such that $\ord_Y D \not = 0$. For each $D$, its support is a finite set.
\item Let $Y_1, \cdots, Y_n$ be Weil divisors on $X$. The \emph{supremum} of $Y_1, \cdots, Y_n$, denoted by $\sup_{1 \leq i \leq n} Y_i$, is a divisor on $X$ defined to be
$$
\sup_{1 \leq i \leq n} Y_i=\sup \{ Y_1, \cdots, Y_n \} := \sum_{Y \in {\rm PDiv} (X)} ( \max_{1 \leq i \leq n } \ord_Y (Y_i) ) [Y].
$$
This expression makes sense for only finitely many $Y \in {\rm PDiv} (X)$ could satisfy $\ord_Y (Y_i) \not = 0$ for some $1 \leq i \leq n$. (\emph{c.f.} Remark \ref{difference})
\end{enumerate}
}
\end{definition}

\subsection{Cubical additive higher Chow complex}
Let $m \geq 2, n \geq 0$ be integers. Let $X$ be an irreducible variety over $k$. Let
\begin{eqnarray*}
\tuborg \square := \bbP^1 \backslash \{ 1 \}, \\ \Diamond_n := \bbA^1 \times \square^n \ni (x, t_1, \cdots, t_n), \\ \widehat{\Diamond}_n := \bbA^1 \times (\bbP^1 )^n \ni (x, t_1, \cdots, t_n ). \sluttuborg
\end{eqnarray*}
For each $i \in \{ 1, \cdots, n\}$ and $j \in \{ 0, \infty\}$ we have codimension one faces $\mu_i ^j : X \times \Diamond_{n-1} \hookrightarrow  X \times \Diamond_n,$ where $(y, x, t_1, \cdots, t_{n-1})\mapsto  (y, x, t_1, \cdots, t_{i-1}, j, t_{i}, \cdots, t_{n-1}).$ 
We also have the ``degenerate'' face $\{x =0 \} \hookrightarrow X \times \Diamond_n$, where $ (y, 0, t_1, \cdots, t_n) \mapsto (y, 0, t_1, \cdots, t_n ).$ Higher codimension faces are obtained by intersecting the above faces.

\begin{definition}\label{cycle groups}Let $F_n \subset X \times \Diamond_n$ be the union of faces $F_i ^j := \mu_i ^j (X \times \Diamond_{n-1})$, $i \in \{ 1, \cdots, n \}$, $j \in \{ 0, \infty \}$.

\begin{enumerate} 
\item [(A)] For $0$-cycles, define
$$c_0 (X \times \Diamond_n) := \bigoplus_{\mathfrak{p}:\mbox{closed pt}} \bbZ \mathfrak{p}, \ \ \ \ \mathfrak{p} \in X \times \Diamond_n \backslash \left( F_n \cup \{x=0\} \right),$$ and for any $m \geq 2$, we let $c_0 (X \times \Diamond_n; m) := c_0 (X \times \Diamond_n)$.

\vskip0.2cm
\item [(B)] For $p$($>0$)-dimensional cycles, define $c_p( X \times \Diamond_n;m )$ inductively as follows.
 
 Suppose that $c_{p-1}(X \times \Diamond_{n-1}; m)$ is defined. Let 
 \begin{eqnarray*} c_p (X \times \Diamond_n;m) := \underset{W}{\bigoplus} ~\bbZ W, \end{eqnarray*} where $W \subset X \times  \Diamond_n$ is an irreducible $p$-dimensional closed subvariety with a normalization $\nu: \overline{W} \to X \times  \widehat{\Diamond}_n$ of the Zariski-closure $\widehat{W}$ of $W$ in $X \times \widehat{\Diamond}_n$, satisfying the following properties:
 
 \begin{enumerate}
 \item [(1)] $W$ intersects all lower dimensional faces properly, that is, in the right codimensions.
 \item [(2)] The associated $(p-1)$-cycle of the scheme $W \cap F_i ^j$ lies in the group $c_{p-1} (X \times \Diamond_{n-1};m)$ for all $i \in \{ 1 , \cdots, n \}$ and $j \in \{ 0, \infty\}$. This cycle is denoted by $\partial_i ^j W$.
 \item [(3)] The Weil divisor \begin{eqnarray}\label{modulus_sup_form}\sup_{1 \leq i \leq n} \left( \nu^* \left\{ t_i = 1 \right\} \right)- m \nu^* \left\{ x=0 \right\}\end{eqnarray} on the normal variety $\overline{W}$ is effective.
 \end{enumerate}
 We obtain the boundary map $\partial:= \sum_{i=1} ^n (-1)^i ( \partial_i ^0 - \partial_i ^{\infty} ) : c_p (X \times \Diamond_n;m) \to c_{p-1} (X \times \Diamond_{n-1};m),$ which satisfies $\partial^2 = 0$.
\end{enumerate}
 \end{definition}

\bigskip

Let $d_p (X \times \Diamond_n;m)$ be the subgroup of $ c_p (X \times \Diamond_n;m)$ generated by degenerate cycles on $X \times \Diamond_n$, \emph{i.e.} those obtained by pulling back admissible cycles on $X \times \Diamond_{n-1}$ via various projection maps. Define $
\calZ_p (X \times \Diamond_n; m) := c_p (X\times \Diamond_n;m) / d_p(X \times \Diamond_n ; m).$

The boundary map $\partial$ on $c_{*} (X \times \Diamond_*;m)$ descends onto $\calZ_{*} (X \times \Diamond_*;m)$ by the usual cubical formalism (\cite{B2,L2,T}).

 \begin{remark}\label{premodulus2}The condition \eqref{modulus_sup_form} of the Definition \ref{cycle groups} is equivalent to the following: for each divisor $Y \in \supp ( \nu ^* \{ x = 0 \} )$ on $\overline{W}$, there is an index $i \in \{ 1 , \cdots, n \}$ for which
\begin{equation}\label{rephrasal}\ord_Y \left( \nu^* \left\{ t_i=1 \right\}  - m \nu^* \left\{ x=0 \right\} \right) \geq 0.\end{equation}
\end{remark}

\begin{definition}\label{modulus2}
\begin{enumerate}
\item [(1)] We say that $W$ has the modulus $m$ condition in $t_i$ along $Y$, and we write $(W, Y) \in \calM^m (t_i)$ if the condition \eqref{rephrasal} holds.
\item [(2)] We write $(W, Y) \in \calM^m (t_i) \cap \calM^m (t_j)$ or $(W, Y) \in \calM^m (t_i, t_j)$, if $(W,Y) \in \calM^m (t_i)$ and $(W, Y) \in \calM^m (t_j)$ for two distinct indices $i, j \in \{ 1 , \cdots, n \}$.
\item [(3)] We generalize (1) as follows. Suppose that $W'$ is any normal variety with a birational surjective morphism $\phi: W' \to \widehat{W}$, where $\widehat{W}$ is the Zariski-closure of $W$ in $X \times \widehat{\Diamond}_n$. By the universal property of normalizations, there is a unique birational morphism $\phi' : W' \to \overline{W}$ such that $\phi = \nu \circ \phi'$. After pulling back along $\phi$, the Weil divisor  \begin{eqnarray*} \sup_{1 \leq i \leq n} \left( \phi^* \left\{ t_i = 1 \right\} \right)- m \phi^* \left\{ x=0 \right\}\end{eqnarray*} is also effective by the condition \eqref{modulus_sup_form}. Hence, for each prime Weil divisor $Y' \in \supp ( \phi^* \{ x = 0 \} )$ on $W'$, there is an index $i \in \{ 1, \cdots, n \}$ such that
\begin{eqnarray*}\ord_{Y'} \left( \phi^* \left\{ t_i=1 \right\}  - m \phi^* \left\{ x=0 \right\} \right) \geq 0.\end{eqnarray*} In this case, we say that $W$ satisfies the modulus $m$ condition in $t_i$ along $Y'$ and we write $(W, Y') \in \calM^m (t_i)$.
\item [(4)] For the above $\phi: W' \to \widehat{W}$ and $i \in \{ 1, \cdots, n \}$,  let
\begin{eqnarray*}
S_i (W') := \supp _i (W'):= \left\{ Y' \in \supp (\phi^* \{ x = 0 \} ) | (W, Y') \in \calM^m (t_i) \right\},
\end{eqnarray*}
\begin{eqnarray*}
S_i ' (W'):=\supp _i ' (W'):=S_i (W') - \bigcup_{j \not = i} S_j (W').
\end{eqnarray*}
Note that we have $
\supp (\phi^* \{ x = 0 \}) = \bigcup_{i=1} ^n \supp_i (W').$

\end{enumerate}
\end{definition}
 
 \begin{definition}\label{ACH}For the complex of abelian groups
 \begin{equation*}\cdots \overset{\partial}{\to} \calZ_{p+1} (X \times \Diamond_{n+1};m) \overset{\partial}{\to} \calZ_{p} (X \times \Diamond_{n}; m) \overset{\partial}{\to} \calZ_{p-1} (X \times \Diamond_{n-1};m) \overset{\partial}{\to} \cdots\end{equation*}define the \emph{cubical additive higher Chow group} $ACH _p(X, n;m)$ to be the homology of the complex at $\calZ_{p} (X \times \Diamond_{n}; m)$: \begin{equation*}ACH_p (X, n;m) := \frac{ \ker ( \partial: \calZ_{p} (X \times \Diamond_{n}; m) \to \calZ_{p-1} (X \times \Diamond_{n-1}; m) )}{{\rm im} ( \partial : \calZ_{p+1} (X \times \Diamond_{n+1}; m)  \to \calZ_{p} (X \times \Diamond_{n}; m))}.\end{equation*}We let $ \calZ^{q} (X \times \Diamond_ n ; m)= \calZ_p (X \times \Diamond_ n; m)$ and $ ACH^{q} (X, n;m)=ACH_p (X, n;m)$, where $q=\dim X + n+1 - p$, in terms of the codimension. Note that the space $\Diamond_n$ has dimension $n+1$. When $m=2$, we use the notation $ACH_p (X, n)$ instead of $ACH_p (X, n; 2)$, and similarly for $\calZ_p (X \times \Diamond_ n)$.
\end{definition}

 \begin{remark}\label{difference}For $X = \spec (k)$, our definition of the group of $1$-cycles $\calZ_1 (\Diamond_n;m)$ differs slightly from the definitions in \cite{BE2,R}. In fact our group is smaller than the ones in \emph{ibids.} where the authors' interest lies on $0$-cycles. For $0$-cycles, this difference does not disturb the main theorem $ACH_0 (k, n;m)\simeq \bbW_{m-1} \Omega_{k} ^{n}$ of \emph{ibids.} (see Remark 3.23-(ii) in \cite{R}). However, for higher dimensional cycles, we believe that the group in \emph{ibids.} maybe is too big. Remark \ref{reason} shows a reason.
 \end{remark}

\section{Regulators on additive higher Chow groups}\label{section:3}

In this section we construct the regulators on the additive higher Chow groups of $1$-cycles for $X = \spec (k)$ where ${\rm char} (k) = 0$.

 \subsection{The statement of the theorem}  
Consider $n$ rational absolute K\"ahler differential $(n-1)$-forms $\omega_{l,m}^n$ ($m, n \geq 2$, $1 \leq l \leq n$) on $\widehat{\Diamond}_n$ defined \emph{cyclically} as follows: 
\begin{eqnarray}\label{Thetas}
{\tuborg \omega_{1,m} ^{n}& = \frac{ 1-t_1}{x^{m+1}} \frac{dt_2}{t_2} \wedge \cdots \wedge \frac{dt_n}{t_n}\\
\omega_{l,m} ^n &= \frac{ 1-t_l}{x^{m+1}} \frac{dt_{l+1}}{t_{l+1}} \wedge \cdots \wedge \frac{dt_n}{t_n} \wedge \frac{dt_1}{t_1} \wedge \cdots \wedge \frac{dt_{l-1}}{t_{l-1}}\ \  (1<l <n) \\
\omega_{n, m} ^n & =  \frac{1-t_n}{x^{m+1}} \frac{dt_1}{t_1} \wedge \cdots \wedge \frac{dt_{n-1}}{t_{n-1}} \sluttuborg
} \end{eqnarray} $$ \in \Gamma\left(\widehat{\Diamond}_n, \Omega_{\widehat{\Diamond}_n/ \bbZ} ^{n-1} (\log F_{n})(* \{ x = 0 \}) \right)$$ When it does not cause confusion, we simply write $\omega_l ^n$ or even $\omega_l$ instead of $\omega_{l,m} ^n$.

 Recall from the Definition \ref{modulus2} that for an irreducible curve $C \in \calZ_1 (\Diamond_n;m)$ and the normalization of its closure $\nu: \overline{C} \to \widehat{C}$, each prime Weil divisor $p \in \supp (\nu^* \{ x = 0 \} )$ is a closed point of $\overline{C}$ and each such $p$ satisfies $(C, p) \in \calM ^m (t_l)$ for at least one $l \in \{ 1 , \cdots, n \}$.

 \begin{theorem}\label{regulator1}Let $k$ be a field of characteristic $0$, and let $m,n \geq 2$ be integers. Then, there is a nontrivial  homomorphism $R_{2,m} : \calZ_1 (\Diamond_n;m) \to \Omega_{k/\bbZ} ^{n-2} $ for which the composition $R_{2,m} \circ \partial : \calZ_2 (\Diamond_{n+1}; m) \overset{\partial}{\to} \calZ_1 (\Diamond_n ; m) \overset{R_{2,m}}{\to} \Omega_{k/\mathbb{Z}} ^{n-2}$ is trivial.
 
This map $R_{2,m}$ is defined as follows: for each irreducible curve $C \in \calZ _1 (\Diamond_n;m)$ with the normalization $\nu: \overline{C} \to \widehat{\Diamond}_n$ as in the Definition \ref{cycle groups}, we let
 \begin{eqnarray}\label{the sum}
 R_{2,m} (C):= \sum_{ p \in \supp ( \nu^* \{ x = 0 \} )} R_{2,m} (C;p), \ \ \mbox{where}\end{eqnarray}
 \begin{eqnarray*} &R_{2,m} (C;p):= (-1)^{l-1} \res_p  \nu^* (\omega_l ^n) ~~~~~~~~~~~~~~~~~~\\ &~~~~~~~ \mbox{ if } (C, p) \in \calM^m (t_l) \mbox{ for some } l  \in \{1, \cdots, n \},\end{eqnarray*}and we extend it $\bbZ$-linearly.
 \end{theorem}

The sum in \eqref{the sum} is finite since the set $\supp ( \nu^* \{ x = 0 \} )$ is finite. An issue is well-definedness; there can be more than one index $l\in \{ 1, \cdots, n \}$ for which $(C, p) \in \calM^m (t_l)$.
  
 \begin{lemma}\label{well-defined1}The map $R_{2,m}$ is well-defined on $\calZ_1 (\Diamond_n;m)$.
\end{lemma}

\begin{proof}We prove that if a closed point $p \in \supp ( \nu^* \{ x = 0 \})$ satisfies $(C, p) \in \calM^m (t_{l}) \cap \calM^m(t_{l'})$ for two distinct indices $l, l' \in \{ 1, \cdots, n \}$, then
\begin{equation*} (-1)^{l -1}\res_p \nu^* (\omega_{l} ^n  ) = (-1)^{l' -1} \res_p \nu^* ( \omega_{l'} ^n ).\end{equation*} We prove that both are zero. We may assume that $l = 1, l' = 2$ and $n = 2$.

Let $t$ be a local parameter at $p \in \overline{C}$. To compute residues, it is sufficient to look at the functions $x, t_1, t_2$ near $p$. By the condition \eqref{rephrasal}, locally near $p$ we have
\begin{equation*}(x(t), t_1 (t), t_2(t)) = (t^r u(t), 1+ t^{rm} f(t), 1+ t^{rm} g(t) )\end{equation*}for some $f(t), g(t) \in \calO_{\overline{C}, p}$ and $u(t) \in \calO_{\overline{C}, p} ^{\times}$, where $r = \ord_p (\nu^* \{ x = 0 \} ) \geq 1$ . By a direct computation near $p$ (keep in mind that $m \geq 2$),
\begin{eqnarray*} \nu^* \left( \frac{1-t_1}{x^{m+1}} \frac{dt_2}{t_2} \right)  &=& \frac{ - t^{rm} f(t)}{t^{rm+r} u(t)^{m+1}} \frac{d ( 1+ t^{rm} g(t) )}{1+ t^{rm} g(t)} \\ &=& - \frac{f(t)}{t^{r} u(t)^{m+1}} \frac{ ( rm t^{rm-1} g(t) + t^{rm} g'(t) ) dt}{1 + t^{rm} g(t)} \\
&=& \frac{ - f(t) ( rm t^{r(m-1)-1} g(t) + t^{r(m-1)} g'(t) ) }{u(t)^{m+1} \left( 1 + t^{rm} g(t) \right)} dt \in \calO_{\overline{C}, p}  dt. ~~~~\end{eqnarray*} 

Hence, \begin{eqnarray*}\res_p \nu^* \left( \frac{1-t_1}{x^{m+1}} \frac{dt_2}{t_2} \right)  = \res_{t=0}  \frac{ - f(t) ( rm t^{r(m-1)-1} g(t) + t^{r(m-1)} g'(t) ) }{u(t)^{m+1} \left( 1 + t^{rm} g(t) \right)} dt  =0.\end{eqnarray*}

Similarly, 
\begin{eqnarray*}\res_p \nu^* \left( \frac{1-t_2} {x^{m+1}} \frac{dt_1}{t_1} \right) = \res_{t=0}  \frac{ - g(t) ( rm t^{r(m-1)-1} f(t) + t^{r(m-1)} f'(t) )}{u(t)^{m+1} (1+ t^{rm} f(t))} dt = 0.\end{eqnarray*} This proves the lemma.
\end{proof}
 
 Hence the points of $\supp ( \nu ^*\{ x = 0 \} )$ of \emph{mixed type}, \emph{i.e.} the points $p \in \overline{C}$ at which $(C, p) \in \calM^m (t_{l}) \cap \calM^m (t_{l'})$ for two distinct indices $l$ and $l'$, play no interesting role in the definition of $R_{2,m}$. This justifies the introduction of notations with prime $S_i '(\overline{C})$ in the Definition \ref{modulus2}-(4).

\begin{remark}\label{reason}Observe the following comparison result before the proof. When $n=m=2$ and $k = \bbC$, the regulator $R_{2,2}$ on $\calZ_1 (\Diamond_2; 2) / \partial \calZ_2 (\Diamond_3;2)$ behaves like the regulator map on the $K$-theoretic version defined in the Proposition 2.3 of \cite{BE2}, that we recall here. 

Let $R = \calO_{\bbA^1, \{0\} }$, $\mathfrak{m} = (x)R$, and let $T\calB_2 (k) = K_2 (R, \mathfrak{m}^2)/ \calC$ for some subgroup $\calC$, where we direct the reader to \cite{BE2} for its precise definition. There is a presentation for $K_2 (R, \mathfrak{m}^2)$ by F. Keune (Theorem 15 in \cite{Keune}, \emph{c.f.} Remark 2.2 in \cite{BE2}) in terms of generators $\left< a, b \right>$ for each $  (a,b) \in (R \times \mathfrak{m} ^2) \cup (\mathfrak{m} ^2 \times R)$ with relations
\begin{eqnarray*}& \tuborg \left< a, b \right> = - \left< b, a \right> ; & a \in \mathfrak{m}^2 \\
 \left<a, b \right> + \left< a , c \right> = \left< a , b + c - abc \right> ; &  a \in \mathfrak{m}^2 \mbox{ or } b ,c \in \mathfrak{m} ^2 \\
 \left< a, bc \right> = \left< ab, c \right> + \left< ac, b \right> ; & a \in \mathfrak{m} ^2 \sluttuborg.\end{eqnarray*} Using this presentation, define $\rho: T\calB_2 (k) \to k\simeq \mathfrak{m} ^2 \Omega_R ^1/ \mathfrak{m} ^3 \Omega_R ^1$ by $\rho( \left<a, b \right>) = -adb$ if $a \in \mathfrak{m}^2$ and $\rho(\left<a, b \right>) = bda$ if $b \in \mathfrak{m}^2$. This is well-defined (Prop. 2.3 in \cite{BE2}).
 
Let us see how $R_{2,2}$ is related to the regulator $\rho$. Since $k = \bbC$, locally we can write $a = \log t_1$, $b= \log t_2$, where $t_1, t_2$ are regarded as functions near $p \in \overline{C}$. If $(C, p) \in \calM^2 (t_1)$, then $t_1= 1 + x^2 f$ for some $f \in \calO_{\overline{C}, p}$ so that $a = \log (1 + x^2 f) = - x^2 f + x^4 (\cdots)$ for some $(\cdots)$. Hence, the regulator value at $p$ is
 \begin{eqnarray*}&\res_p  \nu^*( \frac{1- t_1 }{x^3} \frac{d t_2}{t_2} )  = \res_p \nu^* ( ({ 1- e^a})db/{x^3} )  \\
 & = \res_p  \nu^* ( ({ 1 - 1 - a - \frac{a^2}{2!} - \cdots })db /{x^3} )  = \res_p \nu^* ( { -a db}/{x^3} ) \end{eqnarray*} which coincides with the image of $-a db$ in $k$ under the identification $\mathfrak{m} ^2 \Omega ^1 _R/ \mathfrak{m}^3 \Omega^1 _R \overset{\sim}{\to} k,$ where $ \omega \mapsto \res_{x=0} ( \omega/x^3 ).$ 
 
 Similarly if $(C, p) \in \calM^2 (t_2)$, then $-\res_p  \nu^* ( ((1-t_2)/x^3) (dt_1/t_1) )$ is equal to the image of $bda$ in $k$. When $(C,p) \in \calM^2 (t_1) \cap \calM^2 (t_2)$, Lemma \ref{well-defined1} corresponds to the well-definedness of $\rho$. A precise relationship between the group $T\calB_2 (k)$ and  $\calZ_1 (\Diamond_2; 2) / \partial \calZ_2 (\Diamond_3;2)$ is explained in Corollary 3.2 in \cite{P2}.
 \end{remark}

\subsection{Proof of the Main Theorem} For any irreducible surface $W \in \calZ_2 (\Diamond_{n+1};m)$ we prove that the equality $R_{2,m} (\partial W)=0$ holds.

\subsubsection{Idea of the proof} The details of the proof might seem notationally complicated, but the underlying idea is quite simple.

\begin{figure}\begin{center}
\includegraphics[scale=0.55]{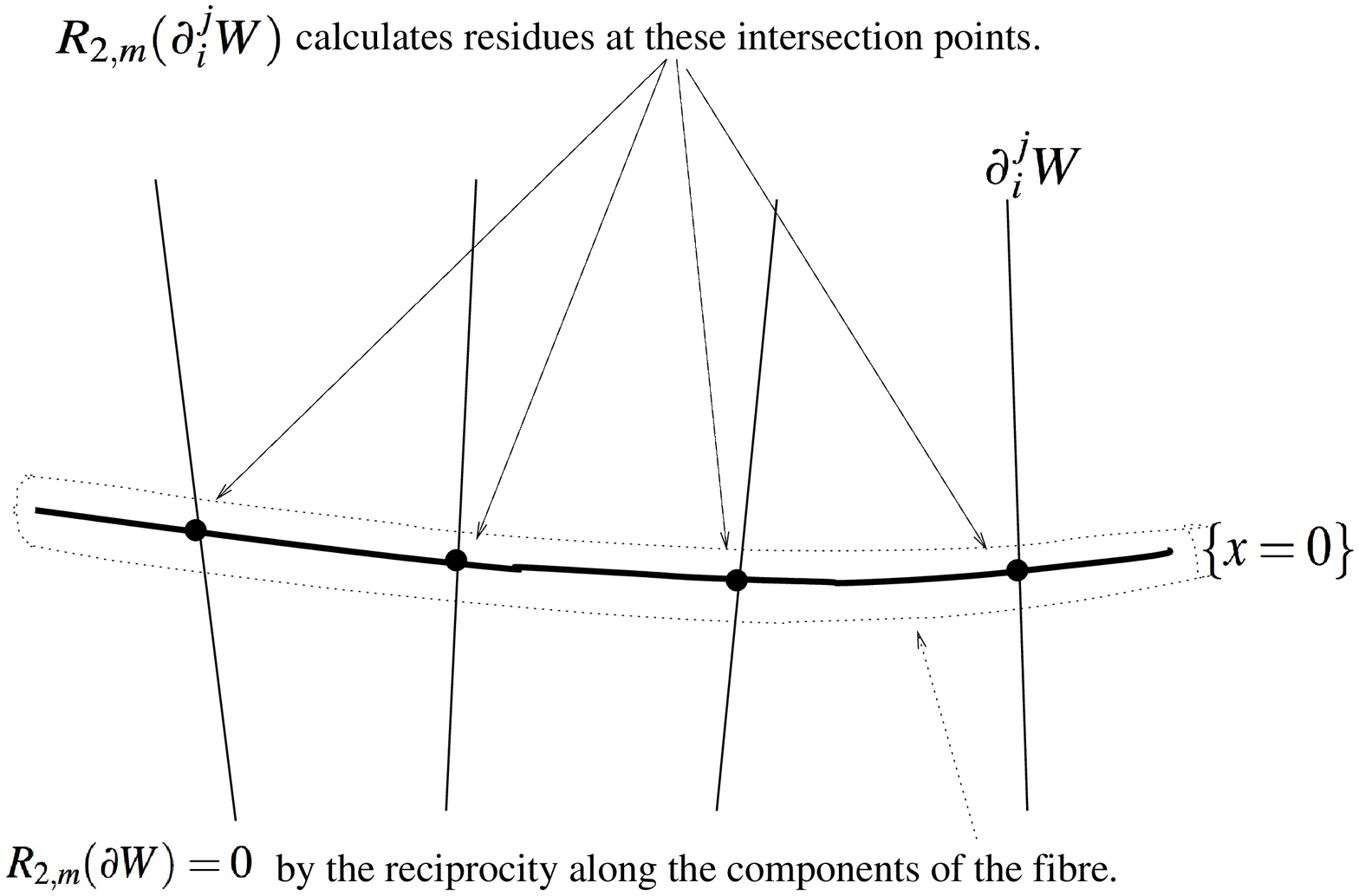}\end{center}
\caption{Sketch of the idea of the proof}\end{figure}

The value $R_{2,m}(\partial W)$ is the alternating sum of all $R_{2,m} (\partial_i ^j W)$ with appropriate signs. Each face $\partial_i ^j W$, depicted as a near-vertical line in the Figure 1, can be seen as cycles on $W$. The regulator value $R_{2,m} (\partial_i ^j W)$ is the sum of residues with some signs of the pull-backs of absolute differential $(n-1)$-forms on $W$ to each component $Y$ of $\partial_i ^j W$, at the points $p$ lying over $\{x=0\}$. These residues of the $(n-1)$-forms are in fact the residues of some absolute differential $n$-forms along the chains $(W, Y, p)$. Since the differential $n$-forms have poles only along the components of $\partial_i ^j W$ and the components over $\{ x=0\}$, using the first part of the reciprocity (Theorem \ref{Parshin-Lomadze}) the value $R_{2,m} (\partial W)$ can be rewritten as the sum of unsigned residues of several differential forms on the union of projective curves over $\{ x =0 \}$. The reciprocity then says that this sum of residues must vanish, proving that $R_{2,m} (\partial W) = 0$.

\subsubsection{Lemmas and the proof}
We use the following version of the resolution of singularity as in \cite{G1,G4}:

\begin{theorem}Let $X$ be a nonsingular quasi-projective variety over a field $k$ of characteristic $0$, and let $V$ be a closed irreducible possibly singular variety in $X$. Let $C$ be an effective cycle on $V$. Then there is a sequence of blow-ups $\widetilde{X} \to \cdots \to X$, whose composition is denoted by $\pi: \widetilde{X} \to X$, such that the strict transform $\widetilde{V}$ of $V$ is nonsingular and $\widetilde{C}:= \pi |_{\widetilde{V}} ^* C$ is a normal crossing divisor whose irreducible components are all nonsingular.
\end{theorem}

\begin{remark}We explain some of our notations. For a variety $Z \subset \Diamond_n$, the Zariski-closure of $Z$ in $\widehat{\Diamond}_{n}$ is denoted by $\widehat{Z}$, the normalization of $\widehat{Z}$ is the morphism $\nu: \overline{Z} \to \widehat{Z}$, and a desingularization of $\widehat{Z}$ is the morphism $\phi:\widetilde{Z} \to \widehat{Z}$.
\end{remark}

For a given irreducible surface $ W \in \calZ_2 (\Diamond_{n+1};m)$, we apply the above desingularization to $X= \widehat{\Diamond}_{n+1}$, $V= \widehat{W}$ and $C = \sum_{ 1 \leq i \leq n+1, j = 0, \infty}  \widehat{\partial_i ^j W}$, where $\widehat{\partial_i ^j W}$ is the cycle on $\widehat{W}$ obtained, firstly, by taking the Zariski-closure of each irreducible component of $\partial_i ^j W$ in $\widehat{\Diamond}_n$, and secondly, by regarding it as an effective $1$-cycle on $\widehat{W}$ via the face map $\mu_i ^j : \widehat{\Diamond}_n \to \widehat{\Diamond}_{n+1}$. Since $\nu: \overline{W} \to \widehat{W}$ is the normalization, by the universal property of normalizations, the desingularization $\phi: \widetilde{W}\to \widehat{W}$ induces a morphism $\phi': \widetilde{W} \to \overline{W}$ such that $\phi = \nu \circ \phi'$.

The morphism $\phi$ gives desingularizations of all irreducible components of the supports of all the faces $\widehat{\partial_i ^j W}$. If necessary, by blowing up some closed points of $\widetilde{X}$, we may assume that each non-exceptional irreducible component of the support of $\widetilde{C}$ is disjoint from all other non-exceptional divisors in the support. Since they are curves, desingularization is equivalent to normalization, which is unique up to isomorphism, so that these blow-ups do not change the regulator value of each component of a face, thus do not change the values of $R_{2,m} (\partial_i ^j W)$. We denote by $\widetilde{\partial_i ^j W}$, the part of the cycle $\widetilde{C}$ whose irreducible components are strict transformations of the irreducible components of the support of $\widehat{\partial_i ^j W}$ via the desingularization $\phi: \widetilde{W} \to \widehat{W}$, that is,
\begin{equation*}
\widetilde{\partial_i ^j W} = \sum _{ \widetilde{Y}: \mbox{ non-exceptional}}(\ord_{\widetilde{Y}} \widetilde{C}) \widetilde{Y}. \end{equation*}

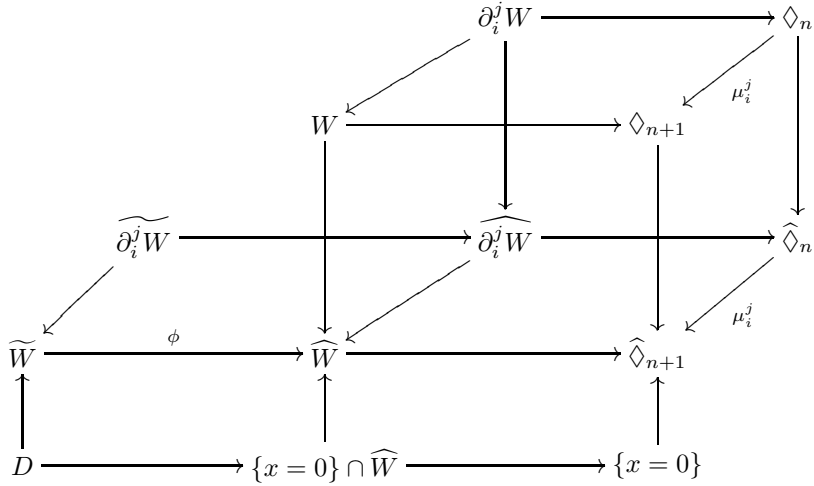
\begin{figure}[htbp]\label{situation1}
$$\xymatrix{ & & &\partial_i ^j W \ar[dl] \ar[dd] \ar[rr] & & \Diamond_n \ar[dl]^{\mu_i ^j} \ar[dd] \\
&&W \ar[dd] \ar[rr] & & \Diamond_{n+1} \ar[dd] & \\
& \widetilde{\partial_i ^j W} \ar[dl] \ar[rr]  & & \widehat{\partial_i ^j W} \ar[dl] \ar[rr] & & \widehat{\Diamond}_n \ar[dl] ^{\mu_i ^j} \\
\widetilde{W} \ar[rr] ^{\phi} & & \widehat{W} \ar[rr] & & \widehat{\Diamond}_{n+1} & \\
D \ar[rr]  \ar[u]  & & \{ x = 0 \} \cap \widehat{W} \ar[u] \ar[rr]& & \{ x = 0 \} \ar[u] &}
$$\caption{After desingularization, where $D \in \supp ( \phi^* \{ x = 0 \} )$}\end{figure}

Let $i \in \{ 1 , \cdots, n+1 \}$ and $j\in \{ 0, \infty \}$. Let $\alpha_1 (i), \cdots, \alpha_n (i)$ be the unique increasing sequence of integers in the set $\{1, \cdots, n+1 \} \backslash \{ i \}$. The embeddings $\mu_i ^j: \widehat{\Diamond}_n \hookrightarrow \widehat{\Diamond}_{n+1}$ maps $t_l$ $(1 \leq l \leq n)$ to $t_{\alpha_l (i)}$, and $\mu_i ^j (\widehat{\Diamond}_n ) = F_i ^j=\{ t _ i = j \}\subset \widehat{\Diamond}_{n+1}$. We can rewrite the definition of the regulator $R_{2,m}$ as follows:

 \begin{corollary}\label{rewritten definition} Under the identification $\mu_i ^j: \widehat{\Diamond}_n \overset{\sim}{\to} F_i ^j=\{ t_i = j \} \subset \widehat{\Diamond}_{n+1}$, the differential $(n-1)$-forms $\omega_l ^n$ $(1 \leq l \leq n)$ on $\widehat{\Diamond}_n$ are the restrictions of the following differential $(n-1)$-forms on $\widehat{\Diamond}_{n+1}$ to $F_i ^j=\{ t_i = j\}:$

\begin{eqnarray*}
  \omega_{1} ^{n} (i) & =& \frac{ 1-t_{\alpha_1 (i)}}{x^{m+1}} \frac{dt_{\alpha_2 (i)}}{t_{\alpha_2 (i)}} \wedge \cdots \wedge \frac{dt_{\alpha_n (i)}}{t_{\alpha_n(i)}} \\
   \omega_{l} ^n (i) &=& \frac{ 1-t_{\alpha_l (i)}}{x^{m+1}} \frac{dt_{\alpha_{l+1} (i)}}{t_{\alpha_{l+1}(i)}} \wedge \cdots \wedge \frac{dt_{\alpha_n(i)}}{t_{\alpha_n (i)}} \wedge \frac{dt_{\alpha_1 (i)}}{t_{\alpha_1 (i)}} \wedge \cdots \wedge \frac{dt_{\alpha_{l-1}(i)}}{t_{\alpha_{l-1}(i)}} \\ & & (1 < l < n ) \\
 \omega_{n} ^n (i)& =&  \frac{1-t_{\alpha_n(i)}}{x^{m+1}} \frac{dt_{\alpha_1 (i)}}{t_{\alpha_1 (i)}} \wedge \cdots \wedge \frac{dt_{\alpha_{n-1}(i)}}{t_{\alpha_{n-1} (i)}} \\
 &\in& \Gamma\left(\widehat{\Diamond}_{n+1}, \Omega_{\widehat{\Diamond}_{n+1}/ \bbZ} ^{n-1} (\log F_{n+1})(* \{ x = 0 \} )\right).\end{eqnarray*}
Furthermore the regulator map $R_{2,m}$ is expressible as follows: for an irreducible curve $C \subset \Diamond_n$ with the normalization $\nu: \overline{C} \to \widehat{C}$, notice that if $(C, p) \in \calM^m (t_l)$, then $(\mu_i ^j (C), p) \in \calM^m (t_{\alpha_l (i)})$. Hence,
 \begin{equation*}
 R_{2,m} (C) = \sum_{p \in \supp (\nu^* \{ x = 0 \})} R_{2,m} (C;p),  \ \ \ \mbox{where}
 \end{equation*}
 \begin{eqnarray*}
 R_{2,m} (C;p)  = R_{2,m} (\mu_i ^j (C); p) = (-1)^{l-1} \res_p ( \nu^* \circ ( \mu_i ^j )^*  \omega_l ^n (i) ) \\ \mbox{   if } ( \mu_i ^j (C), p) \in \calM^m (t_{\alpha_l (i)}).
 \end{eqnarray*}
 \end{corollary}

 \bigskip
 
 Observe that by \eqref{Thetas} we also have $n+1$ rational absolute differential $n$-forms $\omega_l ^{n+1}$ ($1 \leq l \leq n+1$) on $\widehat{\Diamond}_{n+1}$. These $n$-forms are related to the above $(n-1)$-forms $\omega_l ^n (i)$ ($1 \leq l \leq n$) as follows:
 
 \begin{lemma}\label{Theta omega}For indices $l \in \{1, \cdots, n\}$ and $i \in \{ 1, \dots, n+1\}$, we have $
 \omega_{\alpha_l (i)}^{n+1} = (-1)^{i+l} \omega_l ^n (i) \wedge \frac{dt_i}{t_i}.$
 \end{lemma}
 \begin{proof}Up to sign, this is obvious. What matters here is just to keep track of signs. However this is also easy because our differential forms are defined cyclically.
\end{proof}
 
 \begin{lemma}\label{Poincare residue}Let $\sgn(0) := 1$ and $\sgn(\infty):= -1$. Then, the Poincar\'e residue (Remark \ref{order of pole}) along the chain $(\Diamond_{n+1}, F_i ^j)$ of $\omega_{\alpha_l (i)} ^{n+1}$ is
 \begin{equation}\label{Poincare residue identity}
 \res_{(\Diamond_{n+1}, F_i ^j)} ( \omega_{\alpha_l (i)} ^{n+1} ) = - \sgn (j) \cdot (-1)^{i+l} \omega_l ^n (i) |_{F_i ^j}.
 \end{equation}
 \end{lemma}
 \begin{proof}From the shape of the $n$-form $\omega_{\alpha_l (i)} ^{n+1}$, we see that $F_i ^j$ is a simple pole. Hence, it makes sense to talk about Poincar\'e residues. Since $F_i ^j = \{t_i = j \}$ and $\res_{z=0}( dz/z) = 1$, $\res_{z=\infty}( dz/z) = -1$, the equation \eqref{Poincare residue identity} just follows from the Lemma \ref{Theta omega}.
 \end{proof}
 
 \bigskip

 \begin{proof}[\textbf{Proof of Theorem \ref{regulator1}}] Let $W$ be an irreducible surface in $\calZ_2 (\Diamond_{n+1}; m)$. We prove that $R_{2,m} (\partial W) = 0$. Let $R = R_{2,m} (\partial W)$. Let $\phi: \widetilde{W} \to \widehat{W}$ be the desingularization described before. By definition,
 \begin{eqnarray*}
 R &=& \sum_{i=1} ^{n+1} (-1)^i ( R_{2,m} (\partial_i ^0 W) - R_{2,m} (\partial_i ^{\infty} W)) \\ & =& \sum_{i=1} ^{n+1} \sum_{j=0, \infty} \sgn (j) \cdot (-1)^i \cdot R_{2,m} ( \partial_i ^j W) \\
& = & \sum_{i=1} ^{n+1} \sum_{j=0, \infty} \sgn (j) \cdot (-1)^i \sum_{Y \in \supp ( \partial_i ^j W )} (\ord_Y ( \partial_i ^j W) ) R_{2,m} (Y).\end{eqnarray*}

By definition we have $\ord_Y (\partial_i ^j W) = i(Y, W \cap F_i ^j | \Diamond_{n+1})$, the intersection multiplicity of $Y$ (see \cite{F}) in the scheme theoretic intersection $W \cap F_i ^j$ in $\Diamond_{n+1}$. For each $Y \in \supp ( \partial_i ^j W )$, let $\widetilde{Y}$ be its strict transformation via $\phi$. There is an 1-1 correspondence between $\supp ( \partial_i ^j W )$ and $\supp ( \widetilde{\partial_i ^j W} )$ given by $Y \mapsto \widetilde{Y}$, and furthermore $\ord_Y (\partial_i ^j W ) = \ord_{\widetilde{Y}} ( \widetilde{\partial_i ^j W} )$.

So, from the restatement (Corollary \ref{rewritten definition}) of the definition of $R_{2,m}$ and the well-definedness (Lemma \ref{well-defined1}) of $R_{2,m}$, the value $R$ is 
\begin{eqnarray*}
R &=& \sum_{i=1} ^{n+1} \sum_{j=0, \infty} \sum_{\widetilde{Y} \in \supp ( \widetilde{\partial_i ^j W} )}i(Y, W \cap F_i ^j|\Diamond_{n+1})\cdot \sgn (j) \cdot (-1)^i  \\ &\cdot &  \sum_{l=1} ^n \sum_{p \in S ' _{\alpha_l (i)} (\widetilde{Y})} (-1)^{l-1} \res_p ( {\phi |_{\widetilde{Y}}}^* \omega_l ^n (i) )  \end{eqnarray*}\begin{eqnarray*}
& = & \sum_{i=1} ^{n+1} \sum_{j=0, \infty}\sum_{l=1} ^n \sum_{\widetilde{Y} \in \supp ( \widetilde{\partial_i ^j W} )} \sum_{p \in S ' _{\alpha_l (i)}(\widetilde{ Y})} i(Y, W \cap F_i ^j |\Diamond_{n+1})\\ & \cdot &\res_p ( {\phi |_{\widetilde{Y}}}^* ( - \sgn(j) \cdot (-1)^{i+l} \omega_l ^n (i) ) ).
\end{eqnarray*}

Observe that ({\emph{c.f.}} $ \res_{z=0} ( df/f)  = \ord_{z=0} (f)$) $$i(Y, W \cap F_i ^j |\Diamond_{n+1}) \cdot \sgn (j) = \res_{(W, Y)}(\iota^* ( {dt_i}/{t_i}) )= \res_{(\widetilde{W}, \widetilde{Y})} ( \phi^* ( {dt_i}/{t_i} ) ),$$ where $\iota$ is the embedding $W \hookrightarrow \Diamond_{n+1}$, and the second equality is a consequence of Theorem \ref{residue under finite morphism} applied to the desingularization $\phi: \widetilde{W} \to \widehat{W}$.

Hence, by Lemma \ref{Poincare residue} together with the above paragraph, we have
 \begin{equation*}
 R= \sum_{i=1} ^{n+1} \sum_{j=0, \infty} \sum_{l=1} ^n \sum_{\widetilde{Y} \in \supp ( \widetilde{\partial_i ^j W} ) } \sum_{p \in S _{\alpha_l (i)} ' (\widetilde{Y})} \res_p ( \res_{(\widetilde{W}, \widetilde{Y})} \phi^* \omega_{\alpha_l (i)} ^{n+1} ).
 \end{equation*}
 
By applying the transitivity (Lemma 12 in \cite{L4}) of residues, this is equal to
 
\begin{equation}\label{on faces}R = \sum_{i=1} ^{n+1} \sum_{j=0, \infty} \sum_{l=1} ^n \sum_{\widetilde{Y} \in \supp ( \widetilde{\partial_i ^j W} ) } \sum_{p \in S _{\alpha_l (i)} ' (\widetilde{Y})} \res_{(\widetilde{W}, \widetilde{Y}, p)} ( \phi^* \omega _{\alpha_l (i)} ^{n+1} ),\end{equation}
 which is a sum of residues of $\omega_{i'} ^{n+1}$ ($1 \leq i' \leq n+1$) at points lying over $|\phi^* \{ x = 0 \} |$. 

From the shapes of $\omega_{i'} ^{n+1}$ (see \eqref{Thetas}), we notice that $\phi^* \omega_{i'} ^{n+1}$ can possibly have poles only along $D \in \supp ( \phi^* \{ x = 0 \} )$ or along the components of the faces $\widetilde{Y} \in \supp (\widetilde{\partial_i ^j W})$, $i \in \{ 1, \cdots, n+1 \}$, $j \in \{ 0, \infty \}$. Thus, out of all saturated chains of length $2$ on $\widetilde{W}$ ending with a closed point on $|\phi^* \{ x = 0 \} |$, just for $\xi= (\widetilde{W}, D, p)$ or $(\widetilde{W}, \widetilde{Y}, p)$, can we have possibly nonzero residue values of $\res_{\xi} (\phi^* \omega_{i'} ^{n+1})$. Furthermore, the closed points $p$ are necessarily lying in an intersection of some $\widetilde{Y}$ and some $D$. In the equation \eqref{on faces}, notice that for each $i' \in \{ 1, \cdots, n+1 \}$, we are adding up all residues of the form $\phi^* \omega_{i'} ^{n+1}$ along all such chains of the form $(\widetilde{W}, \widetilde{Y}, p)$. Hence, by the Theorem \ref{Parshin-Lomadze}-(1), we can rewrite \eqref{on faces} according to the value $i' = \alpha_l (i) \in \{ 1,\cdots, n+1 \}$ to obtain

 \begin{equation}\label{on x=0}
 R=- \sum_{i'=1} ^{n+1} \sum_{D \in S _{i'}  (\widetilde{W})} \sum_{p \in D} \res_{(\widetilde{W}, D, p)} ( \phi^* \omega_{i'} ^{n+1} ).
 \end{equation}

Now, choose a coefficient field $\sigma_D$ (a theorem of I. S. Cohen, \emph{c.f.} Thm. II-8.25A in \cite{H}) for each $D$, and define $ \Psi ^{i'} _{{D}, \sigma_{D}}:= \res_{(\widetilde{W}, D), \sigma_D} ( \phi^* \omega_{i'} ^{n+1} ).$

By the transitivity of residues (Lemma 12 in \cite{L4}) again, the summand of \eqref{on x=0} is,

\begin{eqnarray*} \res_{(\widetilde{W}, D, p)} ( \phi^* \omega_{i'} ^{n+1}) & =&  ( \res_{(D, p)} \circ \res_{(\widetilde{W}, D), \sigma_{D}} ) ( \phi^* \omega_{i'} ^{n+1} ) \\ &= &\res_{(D, p)} ( \Psi^{i'}_{{D}, \sigma_{D}} ) = \res_p ( \Psi^{i'} _{{D}, \sigma_{D}} ).\end{eqnarray*}

Thus, the value $R$ in the equation \eqref{on x=0} is

\begin{equation*}
R = - \sum_{i'=1} ^{n+1} \sum_{D \in S_{i'} (\widetilde{W})} \sum_{p \in D} \res_{p} ( \Psi ^{i'} _{D, \sigma_{D}} ),
\end{equation*}
but, the last sum $\sum_{p \in D} \res_p ( \Psi ^{i'} _{D, \sigma_{D}})$ is $0$ by applying the Theorem \ref{Parshin-Lomadze}-(2) to the projective but possibly singular curve $D$. Thus $R=0$.\end{proof}

\subsection{Example and Remarks}
\begin{example}Let's compute an example and test the validity of the theorem $R_{2,m} (\partial W) = 0$. For simplicity, take $n=m=2$. This is actually the most interesting case. Consider the following parametrized surface in $\calZ_2 (\Diamond_3; 2)$ which is easily seen to be Zariski-closed:
\begin{equation*}\Sigma = \left\{\left( x, t, \frac{ (1-a_1 x)(1-a_2 x)}{1- (a_1 + a_2)x} , \frac{x+t+2}{x+t+1} \right) \right\}\subset \Diamond_3, \end{equation*} where $x, t$ are parameters and $a_1, a_2 \in k \backslash \{ 0, 1 \}$ with $a _1 + a_2 \not = 0$. This is a rational surface satisfying the modulus $2$ condition in $t_2$. Its faces are computed explicitly as follows, where either $x$ or $t$ is the parameter for each curve.

$$\partial_1 ^0 \Sigma = \left\{\left( x, \frac{ (1-a_1 x)(1-a_2 x)}{1- (a_1 + a_2 )x}, \frac{ x+2}{x+1} \right)\right\},$$ $$
\partial _1 ^{\infty}\Sigma = \left\{\left( x , \frac{ (1- a_1 x)(1-a_2 x)}{1- (a_1 + a_2 )x}, 1 \right)\right\}=0,$$ $$
 \partial _2 ^0 \Sigma= \left\{\left( \frac{1}{a_1}, t, \frac{ \frac{1}{a_1} + t + 2}{\frac{1}{a_1} + t + 1 } \right)\right\} + \left\{\left( \frac{1}{a_2} , t , \frac{ \frac{1}{a_2} + t + 2 }{\frac{1}{a_2} + t + 1} \right) \right\}, $$ $$
\partial _2 ^{\infty} \Sigma = \left\{\left( \frac{1}{a_1 + a_2 } , t , \frac{ \frac{1}{a_1 + a_2} + t + 2 }{ \frac{1}{a_1 + a_2} + t + 1 } \right)\right\},$$ $$
\partial _3 ^0 \Sigma= \left\{\left( x , -x -2 , \frac{ (1-a_1 x )(1-a_2 x)}{1- (a_1 + a_2 )x} \right)\right\}, $$ $$
\partial _3 ^{\infty} \Sigma =\left\{ \left( x, -x -1, \frac{ (1- a_1 x )(1-a_2 x)}{1- (a_1 + a_2 )x} \right)\right\}.$$

We compute $R_{2,2} (\partial_i ^j \Sigma)$ for all $i\in \{ 1,2,3\}$ and $j\in \{ 0, \infty \}$. Obviously, we have $ R_{2,2} (\partial_1 ^{\infty}\Sigma ) = R_{2,2} (\partial _2 ^0 \Sigma) = R_{2,2} (\partial _2 ^{\infty}\Sigma ) = 0.$

The cycle $\partial_1 ^0 \Sigma$ has the modulus $2$ condition in $t_1$, so we use the pull-back of the form $((1-t_1)/x^3)(dt_2 / t_2)$ which is $$- \frac{a_1 a_2} {x (1- (a_1 + a_2)x)} \left( \frac{dx}{x+2} - \frac{dx}{x+1} \right).$$ Thus the residue at $x=0 $ is $R_{2,2} (\partial_1 ^0 \Sigma) =- a_1 a_2 \left( (1/2)  - 1 \right) = a_1 a_2/2$. Similarly, on $\partial _3 ^0 \Sigma$ we use $$- \frac{ - a_1 a_2}{x (1- (a_1 + a_2 )x)} \frac{dx}{x+2},$$ and on $\partial _3 ^{\infty} \Sigma$  we use $$- \frac{ - a_1 a_2}{x (1- (a_1 + a_2)x)} \frac{dx}{x+1}$$ so that $R_{2,2} \left( \partial_3 ^0 \Sigma - \partial_3 ^{\infty} \Sigma \right) = - a_1 a_2/2$. Hence, we have $R_{2,2} (\partial \Sigma) =  (a_1 a_2/2) - (a_1 a_2/2)=0$ as we proved.

In particular, that $R_{2,2} (\partial_1 ^0 \Sigma) = a_1 a_2/2 \not = 0$ shows that $R_{2,2}$ is a nontrivial homomorphism. \end{example}

\begin{remark}It still remains to see how this cubical version is related to the simplicial version in \cite{BE1}. Both are examples of Chow groups for degenerate configurations. The question of how to define such a thing in general is open. The method of \cite{B2} for the usual higher Chow groups that compares the cubical and the simplicial versions depends on the $\bbA^1$-homotopy of the functors $CH_*(\cdot, *)$ and it can't be applied here.
\end{remark}

\begin{question}\label{question1} Are the regulator maps $ R_{2,2} : ACH_1 (k, n)\to \Omega_{k/\bbZ} ^{n-2}$ isomorphisms?
\end{question}

The result \eqref{Hesselholt} suggests that the answer seems very affirmative. It also seems that there exist \emph{higher regulators} $R_{p+1,2}: ACH_{p} (k, n) \to \Omega_{k/ \bbZ} ^{n-2p},$ and we wonder if they are isomorphisms as well. Notice that the right hand side vanishes if $n<2p$, but it is not clear if so does the left hand side. This conjectural identification indicates that we may have the following analogue of the Beilinson-Soul\'e vanishing conjecture for the additive higher Chow groups: 
\begin{conj}Let $X= \spec (k)$ and $\rm{char} (k)=0$. If $n<2p$, then $ACH_{p} (k, n) = 0.$ Equivalently,  if $n > 2j-2$, then $ACH^j (k, n) = 0.$
\end{conj} The author believes that the additive Chow groups are related to the Hochschild homology. For $m \geq 2$, the regulator map $R_{2,m}$ constructed here may lift to $\widetilde{R}_{2,m}: ACH_1 (k, n; m) \to \bbW_{m-1} \Omega_k ^{n-2},$ provided a more general theory of residues with values in the de Rham-Witt forms, which not yet exists.

\bibliographystyle{amsplain}

\end{document}